\documentclass[draft]{amsart}

\usepackage{graphics}
\usepackage[all]{xy}

\usepackage{amssymb,amsmath}
\usepackage{psfrag}

\usepackage{amsfonts}
\usepackage{amscd}
\usepackage{mathrsfs}

\usepackage{graphicx}
\usepackage{epstopdf}
\usepackage{tikz}
\usetikzlibrary{arrows,automata}

\usepackage{graphicx,epsfig,float}

\usepackage{url}

\makeatletter
\def\url@leostyle{%
  \@ifundefined{selectfont}{\def\UrlFont{\sf}}{\def\UrlFont{\small\ttfamily}}}
\makeatother
\newtheorem{thm}{Theorem}[section]
\newtheorem{lem}[thm]{Lemma}

\newtheorem{cor}[thm]{Corollary}
\newtheorem{fact}[thm]{Fact}

\newtheorem{clm}[thm]{Claim}

\newtheorem{crt}[thm]{Criterion}

\newtheorem{defn}[thm]{Definition}

\newtheorem{nrmk}[thm]{Remark}

\newcommand{\pf}{{\bf Proof. }}

\renewcommand{\tilde}{\widetilde}
\renewcommand{\bar}{\overline}

\newcommand{\ZZ}{\mathbb{Z}}

\newcommand{\A}{\mbox{$\mathcal{A}$}}

\newcommand{\Pp}{\mbox{$\mathcal{P}$}}

\newcommand{\bJ}{\mbox{$\mathbf J$}}
\newcommand{\mJ}{\mbox{$\mathbb J$}}

\renewcommand{\mod}{\mathrm{Mod}}

\newcommand{\cov}{\mathrm{Cov}}

\newcommand{\op}{\mathrm{Op}}

\newcommand{\df}{\mathrm{def}}
\newcommand{\Df}{\mathrm{Def}}
\newcommand{\Dfs}{\mathrm{Def}({\mathbb S})}
\newcommand{\tDf}{\tilde{\mathrm{Def}}}

\newcommand{\Ho}{\mathrm{Hom}}

\newcommand{\id}{\mathrm{id}}

\newcommand{\imin}[1]{#1^{-1}}
\newcommand{\lind}[1]{\underset{#1}{\underrightarrow{\lim}}}
\newcommand{\Lind}{\underrightarrow{\lim}}  

\newcommand{\exs}[3]{0 \to {#1} \to {#2} \to {#3} \to 0}

\begin{document}

\title {Invariance of  o-minimal  cohomology with definably compact supports}

\author {M\'{a}rio J. Edmundo}

\address{ Universidade Aberta\\ Rua Braamcamp 90\\
   1250-052 Lisboa, Portugal\\
   and\\
    CMAF Universidade de Lisboa\\
Av. Prof. Gama Pinto 2\\
1649-003 Lisboa, Portugal}

\email{mjedmundo@fc.ul.pt}

\author{Luca Prelli}
\address{ CMAF Universidade de Lisboa\\
Av. Prof. Gama Pinto 2\\
1649-003 Lisboa, Portugal}


\email{lmprelli@fc.ul.pt}

\date{\today}
\thanks{The first author was supported by Funda\c{c}\~ao para a Ci\^encia e a Tecnologia, Financiamento Base 2008 - ISFL/1/209.  The second author was supported by Marie Curie grant PIEF-GA-2010-272021. This work is part of the FCT project PTDC/MAT/101740/2008 and PTDC/MAT/122844/2010.\newline
 }


\begin{abstract}
In this paper we find  general criteria for invariance      
and finiteness results for o-minimal cohomology in an arbitrary o-minimal structure. We apply our criteria and obtain new invariance      
and finiteness  results for  o-minimal cohomology in o-minimal expansions of ordered groups and for  the o-minimal cohomology of  definably compact definable groups in arbitrary o-minimal structures. 
\end{abstract}

\maketitle

\begin{section}{Introduction}\label{section intro}
In this paper we find  general criterion for invariance  (Criterion  \ref{crt inv without c})  
and finiteness (Criterion \ref{crt wilder tool}) results for o-minimal cohomology in an arbitrary o-minimal structure. We apply our criteria and obtain new invariance      
and finiteness  results for  o-minimal cohomology in o-minimal expansions of ordered groups and for the o-minimal cohomology of  definably compact definable groups in arbitrary o-minimal structures. 

In o-minimal expansions of ordered groups, our invariance results (Corollaries \ref{cor inv and finite cp exp groups} (1), \ref{cor inv field without c}, \ref{cor inv  exp groups with c} and \ref{cor inv fields}) extend results known previously only in special cases, namely: (1) Delfs invariance       results in real closed fields (\cite[Theorem 6.10]{D3}) for semi-algebraic cohomology with coefficients in sheaves (relative to the semi-algebraic site), both without supports and with semi-algebraically complete supports; (2) invariance       
results  for o-minimal cohomology, without supports,  with constant coefficients on definable sets (i.e. affine definable spaces) equipped with the o-minimal site  in o-minimal expansions of real closed fields (\cite{bo2} and \cite{ew});  (3) invariance       results  for o-minimal cohomology, without supports,  with constant coefficients on closed and bounded definable sets (i.e. affine definably compact definable spaces) equipped with the o-minimal site in o-minimal expansions of ordered groups (\cite{bf}). 

One should  note that in general, in arbitrary o-minimal expansions of ordered groups,  definable spaces (even definably normal, definably compact ones) need not be affine as in o-minimal expansions of fields. See \cite{el} and \cite[Chapter 10, (1.8)]{vdd}.\\

Another useful case where our invariance   criterion applies is the following:

\begin{thm}\label{thm inv  def groups}
Suppose that ${\mathbb M}$ is an arbitrary  o-minimal structure. Let $G$ be a definably compact definable group. Let $F$ be a sheaf  on the o-minimal site  on $G$. If ${\mathbb S}$ is an elementary extension of ${\mathbb M}$ or  an o-minimal expansion  of ${\mathbb M},$ then
we have
$$H^*(G;F)\simeq H^*(G({\mathbb S});F({\mathbb S})).$$
\end{thm}

Theorem \ref{thm inv  def groups} is an important step towards the computation of the o-minimal cohomology of definably compact definable groups in arbitrary o-minimal structures which are expected to be similar to the  o-minimal cohomology of  definably compact definable groups definable in o-minimal expansions of fields (\cite{eo}).

Similarly, in o-minimal expansions of ordered groups, our finiteness result (Corollaries \ref{cor inv and finite cp exp groups} (2)) for the o-minimal cohomology (with coefficients in a finitely generated module over a noetherian ring)  of a (Hausdorff) definably compact definable space extends  a result known previously only in a special case, namely,   for closed and bounded definable sets (i.e. affine, necessarily definably normal,  definably compact definable spaces) in o-minimal expansions of ordered groups. See \cite{bf} and also \cite{bo2} in the case of   o-minimal expansions of fields.\\

Here we also obtain:

\begin{thm}\label{thm omin gp wilder}
Suppose that ${\mathbb M}$ is an arbitrary o-minimal structure. Suppose that $L$ is a finitely generated module over  a noetherian ring. Let $G$ be a definably compact definable group. Then, for each $p$, $H^p(G;L_G)$ is  finitely generated. 
\end{thm}

This theorem together with Theorem \ref{thm inv  def groups}  is another crucial step  towards the computation of the o-minimal cohomology of definably compact definable groups in arbitrary o-minimal structures which as mentioned above are expected to be similar to the  o-minimal cohomology of  definably compact definable groups definable in o-minimal expansions of fields (\cite{eo}).\\


Let us now compare the ideas involved in the proofs  in this paper with those present in the papers cited in the bibliography. We obtain our  invariance results for o-minimal sheaf cohomology with definably compact supports from our Criterion  \ref{crt inv without c} for invariance of o-minimal sheaf cohomology without  supports in the same way Delfs obtained his result \cite[Theorem 6.10 (i)]{D3} for invariance of semi-algebraic  sheaf cohomology with semi-algebraically complete supports from his comparison result \cite[Theorem 6.1]{D3} for  semi-algebraic sheaf cohomology with supports on an arbitrary  semi-algebraic family of supports. Regarding the results used in this reduction, we use the o-minimal analogues, proved in \cite{ep1},  of those used in the semi-algebraic case. Delfs comparison result \cite[Theorem 6.1]{D3} is much stronger than our Criterion  \ref{crt inv without c} and to prove it Delfs uses, in the second part of the proof (\cite[Proposition 6.4]{D3}), the semi-algebraic triangulation theorem in a very involved way. The same method applies as well in o-minimal expansions of fields using the o-minimal triangulation theorem  (\cite[Chapter 8, (2.9)]{vdd}) instead, but, of course, this method does not generalize to the general case covered here. In the first part of Delfs proof, \cite[Proposition 6.3]{D3}, there is a general method used to reduce a result on semi-algebraic sheaves to the case of constant semi-algebraic sheaves on affine semi-algebraic spaces and  there is a use of the semi-algebraic triangulation theorem to prove the result in the later case.  In the proof of our Criterion \ref{crt inv without c}  we borrow this method, see Fact \ref{fact criteria for o-min sheaves}, but of course, since our Criterion \ref{crt inv without c} is not the o-minimal analogue of \cite[Proposition 6.3]{D3}, we have to make some extra but very standard homological algebra computations  before our proof, namely the lemmas before  Criterion \ref{crt inv without c}.

The basic assumption of Criterion \ref{crt inv without c} is the invariance of o-minimal cohomology with constant coefficients of affine definable spaces. 
So in the applications we only need to verify these basic assumptions. For our applications to invariance results in o-minimal expansions of fields and in o-minimal expansions  of ordered groups these basic assumptions were proved in the o-minimal literature as we mentioned before. For  the invariance result for definably compact definable groups in arbitrary o-minimal structures, Theorem \ref{thm  inv  def groups}, we use one of the main results of \cite{epr} 
to reduce its proof to  proving the invariance of o-minimal cohomology with constant coefficients of affine definably compact definable spaces in cartesian products of definable group-intervals. For the proof of the later we adapt in a straight forward  way the  proof in \cite{bf} of the invariance of o-minimal cohomology with constant coefficients of affine definably compact definable spaces in o-minimal expansions of ordered groups.\\

Regarding our finiteness results for o-minimal cohomology of definably compact definable spaces with coefficients in finitely generated modules over noetherian rings, our Criterion \ref{crt wilder tool}, which reduces the problem to affine definably compact definable spaces, is obtained adapting the method used in the proof of Wilder's finiteness theorem (\cite[III.10]{i}). Our applications of this Criterion are obtained as above, using the affine version already present in the o-minimal literature (\cite{bf} and also \cite{bo2}),  or adapting the proof of \cite{bf} to the setting of cartesian products of definable group-intervals.\\

\medskip
\emph{Acknowledgements.} We wish to thank the referee for his/her patient and very helpful work: he/she suggested simplifications of definitions and many proofs as well as considerable improvements in the presentation of  the paper. 

\end{section}

\begin{section}{Preliminaries}\label{section prelim}
In this Section  we will recall  the notions that will be used later. We let ${\mathbb M}$ be an arbitrary o-minimal structure and definable means definable in ${\mathbb M}$ possibly with parameters. \\

We will work in paper in the category $\Df$ of definable spaces with continuous definable maps. See \cite[Chapter 10,  \S 1]{vdd} for the definition and other basic properties and notions about definable spaces.\\

A  definable space $X$  is \textit{definably normal} if one of the following equivalent conditions holds:
\begin{enumerate}
\item
for every disjoint closed definable subsets $Z_1$ and $Z_2$ of $X$ there are disjoint open definable subsets $U_1$ and $U_2$ of $X$ such that $Z_i\subseteq U_i$ for $i=1,2.$

\item
for every  $S\subseteq X$ closed definable  and $W\subseteq X$ open definable such that $S\subseteq W$, there is an open definable subsets $U$  of $X$ such that $S\subseteq U$ and $\bar{U}\subseteq W$.
\end{enumerate}

As usual we have (compare with \cite[Chapter 6, (3.6)]{vdd}):

\begin{fact}[The shrinking lemma]\label{fact shrinking lemma}
Suppose that $X$ is a definably normal definable  space. If $\{U_i:i=1,\dots ,n\}$ is a covering of $X$ by open definable subsets, then there are definable open subsets $V_i$ and definable closed subsets $C_i$ of $X$ ($1\leq i\leq n$) with $V_i\subseteq C_i\subseteq U_i$ and $X=\cup \{V_i:i=1,\dots, n\}$.\\
\end{fact}

\begin{defn}\label{defn P affine on X}
{\em
Let $X$ be a definable space. We say that a subset $C\subseteq X$ is {\it affine} if as a definable subspace it is definably homeomorphic to a definable subset of some $M^k$ with the induced topology.
}
\end{defn}

In this paper we shall use the following general strategy. To show that a  property $\Pp$ of definable spaces holds on a definably normal definable space $X$ we: (a) use the shrinking lemma to  show that if $\Pp$ holds in every  affine closed definable subset of $X$ then $\Pp$ holds on $X$; (b) show $\Pp$ for affine closed definable subsets of $X$. Of course, in each case the proofs of (a) and (b) are specific to the given $\Pp$. Observe however that not every $\Pp$ that holds for affine closed definable subsets of $X$ holds for $X$, e.g. the property saying the space is affine 
since there are non affine definable spaces. Often we will be able to establish  only (a) in general obtaining a criterion for proving $\Pp$ is some interesting cases where $\Pp$ holds for affine closed definable subsets.\\


Let $X$ be a definable space and $C\subseteq X$ a definable subset. By a {\it definable curve in $C$} we mean a continuous definable map  $\alpha :(a,b)\to C\subseteq X$, where $a<b$ are in $M\cup \{-\infty, +\infty \}.$ We say that a definable curve $\alpha :(a,b)\to C\subseteq X$ in $C$ is {\it completable in C} if both limits  $\lim _{t\to a^+}\alpha (t)$ and $\lim _{t\to b^-}\alpha (t)$ exist in $C.$  We say that $C$ is {\it definably compact} if every definable curve in $C$ is completable in $C.$ See \cite{ps} or \cite[Chapter 6]{vdd} in the affine case in o-minimal expansions of ordered groups.\\

With this definition we have that a definable set $X\subseteq M^n$ with its induced topology is definably compact if and only if it is closed and bounded in $M^n$ (\cite[Theorem 2.1]{ps}).  \\
Using the general strategy mentioned above and \cite[Theorem 2.1]{ps} we easily see that:

\begin{nrmk}\label{nrmk def comp2}
Suppose that $X$ is a definably normal definable space. If $K$ is definably compact subset of $X$, then $K$ is a closed definable subset.\\
\end{nrmk}


Recall also that if ${\mathbb S}$ be an elementary extension of ${\mathbb M}$ or  an o-minimal expansion of ${\mathbb M},$ then as it is well known  ${\mathbb S}$  determines a functor from the category of (${\mathbb M}$-)definable sets and (${\mathbb M}$-)definable maps to the category of ${\mathbb S}$-definable sets and ${\mathbb S}$-definable maps. This functor extends to a functor  $\Df \to \Df ({\mathbb S})$ sending a definable space $X$ to the ${\mathbb S}$-definable space $X({\mathbb S})$ and sending a definable (resp. continuous definable) map $f:X\to Y$ to the   ${\mathbb S}$-definable (resp. continuous ${\mathbb S}$-definable) map $f^{{\mathbb S}}:X({\mathbb S})\to Y({\mathbb S})$. The functor $\Df \to \Df ({\mathbb S})$ is a monomorphism from  the boolean algebra of definable subsets of  a definable space $X$ to  the boolean algebra of ${\mathbb S}$-definable  subsets of $X({\mathbb S})$ and it commutes with: (i) the interior and closure operations; (ii) the  image and inverse image under (continuous) definable maps.\\


Since for a given  definable family of definable curves the existence of limits is a first-order condition on the parameters of the family we have:

\begin{nrmk}\label{nrmk def normal inv of compact1}
Let ${\mathbb S}$ be an elementary extension of ${\mathbb M}$. Let $X$ be a definably normal, definably compact definable space. Then $X({\mathbb S})$ is an  ${\mathbb S}$-definably compact ${\mathbb S}$-definable space.
\end{nrmk}

Using the general strategy mentioned before and \cite[Theorem 2.1]{ps} 
we also have:

\begin{nrmk}\label{nrmk def normal inv of compact2}
Let ${\mathbb S}$ be  an o-minimal expansion of ${\mathbb M}$.  Let $X$ be a definably normal, definably compact definable space. Then $X({\mathbb S})$ is an  ${\mathbb S}$-definably compact ${\mathbb S}$-definable space.\\
\end{nrmk}

Finally we observe that later in the paper we will also work occasionally with the category $\tDf$ whose objects are o-minimal spectra of definable spaces and whose morphisms are the corresponding o-minimal spectra  of  (continuous) definable maps. Recall that: (i) the o-minimal spectrum $\tilde{X}$ of a definable space $X$ is, as in the affine case (\cite{c}, \cite{cr} and \cite{p}),  the set of ultrafilters  of definable subsets of $X$ (also called in model theory, types concentrated on $X$) equipped with the topology generated by the open subsets of the form $\tilde{U}$, where $U$ is an open definable subset of $X$; (ii) the o-minimal spectrum $\tilde{f}:\tilde{X}\to \tilde{Y}$ of a (continuous) definable map $f:X\to Y$ between definable spaces is the (continuous) map such that given an ultrafilter $\alpha \in \tilde{X}$, $f(\alpha )$ is the ultrafilter in $\tilde{Y}$ determined by the collection $\{A: f^{-1}(A)\in \alpha \}.$

We will refer the reader to \cite{ejp} for   basic results and notions 
about  o-minimal spectral spaces or about the tilde functor $\Df \to \tDf.$ Note that these results were stated in \cite{ejp} in the category of  definable sets but are true in the category of definable spaces with exactly the same proofs. In fact most  of them hold also in real algebraic spaces (\cite{BCR}, \cite{cr}) and more generally in  spectral topological  space (\cite{cc}).\\

\end{section}

\begin{section}{Criteria for invariance      
and finiteness results}\label{section criteria}
In this section we recall basic facts about sheaves on topological spaces, about o-minimal sheaves and we show the  criteria for our invariance      
and finiteness results.\\

\begin{subsection}{Sheaves}
Let $X$ be a topological space, let $\op(X)$ be the category of open subsets of $X$ (morphisms are given by the inclusions) and let $A$ be a  ring. We denote by $\mod(A_X)$ the category of sheaves $A$-modules  on $X$. We will call the objects of $\mod (A_X)$ $A$-sheaves on $X$. Here we recall some notions and some useful facts about $A$-sheaves on topological spaces. These  general notions and results  apply also to $A$-sheaves on objects of $\tDf$. We refer to \cite{b}, \cite{g}, \cite{i} and \cite{ks1} for further details on these results and other results on $A$-sheaves on topological spaces that we will use later for $A$-sheaves on objects of $\tDf$.\\

An $A$-sheaf  on $X$ is a contravariant functor $F:\op(X)^{\rm op} \to \mod(A_X)$, $U \mapsto \Gamma(U;F)$ satisfying following gluing conditions, which are described, for each $U \in \op(X)$ and each covering $\{U_i\}_{i\in I}\subseteq \op(U)$ of $U,$ by:
\begin{itemize}
\item[(S1)]
if $s\in \Gamma (U;F)$ and $s_{|U_i}=0$ for all $i\in I,$ then $s=0.$
\item[(S2)]
if $s_i\in \Gamma (U_i;F),$ $i\in I$ is a family such that $s_{i|U_i\cap U_j}=s_{j|U_i\cap U_j}$ for all $(i,j)\in I^2,$ then there exists $s\in \Gamma (U;F)$ such that $s_i=s_{|U_i}$ for all $i\in I.$
\end{itemize}

A fiber $F_x$ of $F$ on a point $x \in X$ is given by the limit $\lind {x \in U \in \op(X)} \Gamma(U;F)$. A sequence $0 \to F \to G \to H \to 0$ is exact on $\mod(A_X)$ if it is exact on fibers, i.e. if $0 \to F_x \to G_x \to H_x \to 0$ is exact for each $x \in X$.\\


Let $f:X \to Y$ be a continuous map. The direct image  functor $f_*:\mod(A_X) \to \mod(A_Y)$  is defined by $\Gamma(U;f_*F)=\Gamma(f^{-1}(U);F)$ for $F \in \mod(A_X)$. The inverse image functor $f^{-1}:\mod(A_Y) \to \mod(A_X)$ is defined as follows: if $G \in \mod(A_Y)$, then $f^{-1}G$ is the sheaf associated to the presheaf $U \mapsto \lind {U \subseteq f^{-1}(V)} \Gamma(V;G)$. The functor $f^{-1}$ is left adjoint to the functor $f_*$, i.e. we have a functorial isomorphism
$\Ho(f^{-1}G,F)\simeq  \Ho(G,f_*F);$
the direct image functor $f_*$ is left exact and commutes with  small projective limits; the inverse image functor $f^{-1}$ is exact and commutes with  small inductive limits. When $i_U$ is the inclusion of an open subset on $X$ we have $i_U^{-1}F=F_{|U}$.\\

Let $i_Z:Z\to X$ be the inclusion of a locally closed subset $Z$ of $X$. We recall the definition of the functor $i_{Z!}$ ({\it extension by zero}). This functor is such that, for  $F \in {\rm Mod}(A_Z)$, $i_{Z!}F$ is the unique $A$-sheaf in ${\rm Mod}(A_X)$ inducing $F$ on $Z$ and zero on $X\setminus Z$. First let $U$ be an open subset of $X$ and let $F \in {\rm Mod}(A_U)$. Then $i_{U!}F$ is the sheaf associated to the presheaf $V\mapsto \Gamma(V;i_{U!}F)$ which is $\Gamma(V; F)$ if $V \subseteq U$ and $0$ otherwise. If $S$ is a closed subset of $X$ and $F \in {\rm Mod}(A_S)$, then $i_{S!}F=i_{S*}F$. Now let $Z=U \cap S$ be a locally closed subset of $X$, then one defines $i_{Z!}=i_{U!} \circ i_{S!} \simeq i_{S!} \circ i_{U!}$.

If $f:X \to Y$ is a continuous map, $Z$  a locally closed subset of $Y,$
\begin{equation}\label{loc closed base change}
\xymatrix{
f^{-1}(Z)  \ar@{^{(}->}[r]^{j}  \ar[d]^{f_|} & X \ar[d]^{f} \\
Z \ar@{^{(}->}[r]^{i} & Y
}
\end{equation}
a commutative diagram and $G\in \mod (A_{Z}),$ then $f^{-1}\circ i_!G\simeq  j_!\circ (f_|)^{-1}G.$\\

Let $F \in \mod(A_X)$. One sets $F_Z=i_{Z!}\circ i_Z^{-1}F$. Thus $F_Z$ is characterized by $F_{Z|Z}=F_{|Z}$ and $F_{Z|X\setminus Z}=0$. It is an exact functor. If $Z'$ is another locally closed subset of $X$, then $(F_Z)_{Z'}=F_{Z\cap Z'}.$
Let $L$ be an $A$-module. When $F=L_X$ is the constant sheaf on $X$ of fiber $L$ we just set $L_Z$ instead of $(L_X)_Z$. In particular, for an open subset $U$ of $X$, $A_U$ is the notation for $(A_X)_U.$\\

The functor $(\bullet)_Z$ admits a right adjoint, denoted by $\Gamma_Z$ which is left exact. Let $V \in \op(X)$. When $Z=U \in \op(X)$ we have $\Gamma(V;\Gamma_UF)=\Gamma(U \cap V;F)$. When $Z$ is closed $\Gamma(V;\Gamma_ZF)=\{s \in \Gamma(V;F):\;{\rm supp}\,s \subseteq Z\}$ where ${\rm supp} \, s$ is the complement in $V$ of the union of all open sets $U\subseteq V$ such that $s_{|U}=0.$\\

Let $\Phi $ be a family of supports on $X$ (i.e. a collection of closed subsets of $X$ such that: (i)  $\Phi $ is closed under finite unions and (ii) every closed subset of a member of $\Phi $ is in $\Phi$). Recall that for $F \in \mod(A_X)$, an element $s\in \Gamma (X;F)$ is in $\Gamma _{\Phi }(X;F)$ if and only if ${\rm supp}\,s$
is in $\Phi$, i.e.
$$
\Gamma_\Phi(X;F)=\lind {S \in \Phi}\Gamma(X;\Gamma_SF).
$$



Later in the paper we shall use the right derived versions of many of the above formulas relating the various operations on $A$-sheaves. We will use these derived formulas freely and refer to reader to \cite[Chapter II]{ks1} for details. For instance, the cohomology with supports on $\Phi$ is defined by
$$H^*_{\Phi}(X;F)=R^*\Gamma _{\Phi }(X;F).$$

\end{subsection}

\begin{subsection}{O-minimal sheaves}\label{subsection omin sheaves}
Let $A$ be a  ring. If $X$ is an object of $\Df,$ then the {\it o-minimal site} $X_{\df}$ on  $X$ is the category $\op (X_{\df})$ whose objects are open definable subsets of $X$, the morphisms are the inclusions and the admissible covers $\cov (U)$ of $U\in \op (X_{\df})$ are  covers by open definable subsets with finite subcoverings. We will denote by $\mod(A_{X_{\df}})$ the category of sheaves of $A$-modules on $X$. \\

The tilde functor $\Df \to \tDf$ determines a morphism of sites
$$\nu _X:\tilde{X}\to X_{\df}$$
given by the functor
$$\nu _X^t:\op (X_{\df})\to \op (\tilde{X}):U\mapsto \tilde{U}.$$

\begin{thm}[\cite{ejp}]\label{thm main iso on sheaves}
The inverse image of  $\nu _X:\tilde{X}\to X_{\df}$ determines  an isomorphism of categories
$$\mod (A_{X_{\df }}) \to \mod (A_{\tilde{X}}):F\mapsto \tilde{F},$$
where $\mod (A_{\tilde{X}})$ is the category of $A$-sheaves on the topological space $\tilde{X}$.
\end{thm}

The functors $f_*$ and $\Ho _{A_{X_{\df}}}(\bullet, \bullet )$ commute with the tilde functor by definition. From this one can see that 
$f^{-1}$
commutes by adjunction.\\

By Theorem \ref{thm main iso on sheaves} to develop sheaf theory in $\Df$ is equivalent to developing sheaf theory in $\tDf$. For instance, if $X$ is an object of $\Df$ and if  $\Phi $ is a family of definable supports on $X$ (i.e. a collection of closed definable subsets of $X$ such that: (i)  $\Phi $ is closed under finite unions and (ii) every closed definable subset of a member of $\Phi $ is in $\Phi$), then $\tilde{\Phi}$, the collection of all closed subsets of tildes of members of $\Phi$, is a family of supports on $\tilde{X}$ and we set
$$H_{\Phi }^*(X;F)=H^*_{\tilde{\Phi }}(\tilde{X};\tilde{F}).$$

In the paper \cite{ep1} we used this approach to develop the theory of $\Phi$-supported sheaves, where $\Phi$ is definably normal, namely a family of definable supports supports such that: (1) each element of $\Phi$ is definably normal, (2) for each $S \in \Phi$ and each open definable neighborhood $U$ of $S$ there exists a closed definable neighborhood of $S$ in $U$ which is in $\Phi.$  Below we will use this theory and refer the reader to \cite{ep1} for details.

\begin{nrmk}\label{nrmk k ring}
{\em
Note that in \cite{ep1} we assumed that $A$ is a field, but this is only used there when dealing with the tensor product operation $\bullet \otimes _{A_X}G$ on $A$-sheaves (so that it is always exact). Here we will not require this operation.
}
\end{nrmk}

Also it is often useful to use Theorem \ref{thm main iso on sheaves} to define new operations on o-minimal sheaves. For example, we can extend the usual  definition of the extension by zero operation $i_{U!}: \mod (U_{\df})\to \mod (X_{\df})$ on $A$-sheaves on a site where $U\in \op (X_{\df})$, to the extension by zero operation $i_{Z!}: \mod (Z_{\df})\to \mod (X_{\df})$ on $A$-sheaves on a site where $Z$ is a definable locally closed subset of $X$, by setting

$$\tilde{i_{Z!}F}=\tilde{i}_{\tilde{Z}!}\tilde{F}.$$\\

It is useful to recall here the following general criterion:

\begin{fact}\label{fact criteria for o-min sheaves}
Let  $X$ be an object of $ \Df$  and let $\mathfrak{R}$ be a class of objects of $\mod (A_{X_{\df}}).$ Suppose that $\mathfrak{R}$ satisfies:
    \begin{itemize}
    \item[(i)] for each exact sequence $\exs{F'}{F}{F''}$ with $F' \in \mathfrak{R}$ we have $F \in \mathfrak{R}$ if and only if $F'' \in \mathfrak{R}$;
    \item[(ii)] $\mathfrak{R}$ is stable under filtrant $\Lind$;
    \item[(iii)] $A_V \in \mathfrak{R}$ for any $V \in \op(X_{\df}).$ 
    \end{itemize}
Then $\mathfrak{R}=\mod(A_{X_{\df}}).$ Moreover, 
if $X$ has a finite cover by open definable subsets $\{W_i\}_{i=1}^m$ such that each $\bar{W_i}$ is affine,
then we can replace  $\textrm{(iii)}$ by:
\begin{itemize}
\item[(iii)$'$]
$A_U \in \mathfrak{R}$ for any  $U \in \op(X_{\df})$ such that $U\subseteq W_i$ for some $i.$
\end{itemize}
\end{fact}

This is the o-minimal analogue of a similar criterion in the semi-algebraic case (\cite[Lemma 4.18]{D3}) and is obtained by using the isomorphism $\mod (A_{X_{\df}})\to \mod (A_{\tilde{X}})$  of Theorem \ref{thm main iso on sheaves} and applying the corresponding criterion in the topological case (\cite[Chapter II, 16.12]{b}) (point (iii) is a little bit stronger here) observing that constructible open subsets of $\tilde{X}$ form a filtrant basis for the topology of $\tilde{X}.$

On the other hand, suppose that 
$X$ has a finite cover by open definable subsets $\{W_i\}_{i=1}^m$ such that each $\bar{W_i}$ is affine.
Let $V \in \op(X_{\df})$. 
Then 
$V$ has a finite cover $\{U_i\}_{i=1}^m$  consisting of open  definable subsets of $V$ such that $U_i\subseteq W_i.$ 
By (iii)$'$ each $A_{U_i}\in \mathfrak{R}$, and by  induction on $m$, we see that  $A_V\in \mathfrak{R}$ and so we obtain (iii) which  together (i) and (ii) gives that $\mathfrak{R}=\mod(A_{X_{\df}}).$ So let's show that  $A_V\in \mathfrak{R}$  by induction on $m$. If $m=1$ then $V=U_1$ and the result follows by (iii)$'$. For the inductive step, let $V'=U_1\cup \cdots \cup U_{m-1}$. So we have an exact sequence $0\to A_{V' \cap U_m}\to A_{V'} \oplus A_{U_m}\to A_{V}\to 0$ (\cite[Proposition 2.3.6 (vii)]{ks1}) with $A_{V'\cap U_m}\in \mathfrak{R}$ (by (iii)$'$) and $A_{V'} \oplus A_{U_m}\in \mathfrak{R}$ (by inductive hypothesis, (iii)$'$ and the fact that by (i) $\mathfrak{R}$ is stable under finite sums)
Thus, by (i), we have $A_V\in \mathfrak{R}$ as required. \\

By the shrinking lemma:

\begin{nrmk}\label{nrmk crit normal}
Let $X$ be an object of $\Df.$ If  $X$ is a definably normal, then $X$ has a finite cover by open definable subsets $\{W_i\}_{i=1}^m$ such that each $\bar{W_i}$ is affine. \\ 
\end{nrmk}

\end{subsection}

\begin{subsection}{Criterion for invariance  results}\label{subsection invariance      }
Here we prove our general criterion for invariance       of o-minimal sheaf cohomology  without supports. \\


Below we let $\mathbb{S}$ be an elementary extension of $\mathbb{M}$ or  an o-minimal expansion of ${\mathbb M}$. \\


Recall that given $X$ an object of
$ \Df$, there is a continuous surjective map 
$$r:\widetilde{X(\mathbb{S})} \to \widetilde{X}$$
 defined as follows: for each $\alpha \in \widetilde{X(\mathbb{S})}$, $r(\alpha)=\{A:\; A(\mathbb{S})\in \alpha \}$. If $F \in \mod(A_{X_{\df}})$, then the adjunction morphism ${\rm id} \to Rr_* \circ r^{-1}$ together with the isomorphisms $\mod (A_{X_{\df}})\to \mod (A_{\tilde{X}})$ and $\mod (A_{X({\mathbb S})_{\df}})\to  \mod (A_{\tilde{X({\mathbb S})}})$ of Theorem \ref{thm main iso on sheaves}  define a morphism
\begin{equation}\label{eqn invariance}
R\Gamma(X;F) \simeq R\Gamma(\tilde{X};\tilde{F}) \to R\Gamma(\tilde{X};Rr_*r^{-1}\tilde{F})  \simeq R\Gamma(X(\mathbb{S});F(\mathbb{S}))
\end{equation}
where $F(\mathbb{S})\in \mod (A_{X({\mathbb S})_{\df}})$ is the unique object such that $\tilde{F(\mathbb{S})}=\imin r \tilde{F}$. Above, $Rr_*$ is the right derived functor of the direct image functor $r_*$ and $r^{-1}\tilde{F},$ by definition of inverse image functor,  is the sheaf on $\widetilde{X({\mathbb S})}$ associated to the preasheaf $V\mapsto  \lind {V \subseteq \tilde{U({\mathbb S})}} \tilde{F}(\tilde{U}),$ with $V$ open in $\widetilde{X({\mathbb S})}$ and $\tilde{U}$ open constructible   in $\widetilde{X}.$\\

Below we shall use the  Fact \ref{fact  criteria for o-min sheaves} to get our general criterion for  invariance   results. But first we make a couple of observations.\\

\begin{lem}\label{nrmk inv X to U}
Let $X$ be an object of $\Df$ and let $U\in \op (X_{\df}).$ Then we have isomorphisms 
$$H^*(\bar{U};A_U) \simeq H^*(X;A_U)\,\,\textrm{and}\,\, H^*(\bar{U}(\mathbb{S});A_{U(\mathbb{S})}) \simeq H^*(X(\mathbb{S});A_{U(\mathbb{S})}).$$
\end{lem}

\pf
The exact sequence $0 \to A_U \to A_{\overline{U}} \to A_{\bar{U}\setminus U} \to 0$ (\cite[Proposition 2.3.6 (v)]{ks1}) defines a distinguished triangle $A_U \to A_{\bar{U}} \to A_{\bar{U} \setminus U} \to$ (\cite[Proposition 1.7.5]{ks1}). We have the following morphism of distinguished triangles
$$
\xymatrix{
A_U \ar[r] \ar[d] & A_{\bar{U}} \ar[d] \ar[r] & A_{\bar{U}\setminus U} \ar[d] \ar[r] &\\
i_{\bar{U}*}i_{\bar{U}}^{-1}A_{U} \ar[r] & i_{\bar{U}*}i_{\bar{U}}^{-1}A_{\bar{U}} \ar[r] & i_{\bar{U}*}i_{\bar{U}}^{-1} A_{\bar{U}\setminus U} \ar[r] &
}
$$
where the vertical morphisms are determined by the adjunction $\id \to i_{\bar{U}*}i_{\bar{U}}^{-1}$ and $i_{\bar{U}}:\bar{U}\to X$ is the inclusion (remember that $i_{\bar{U}*} \simeq Ri_{\bar{U}*}$ since $i_{\bar{U}*}$ is exact). Applying the triangulated functor $R\Gamma(X;\bullet)$ we obtain the following morphism of distinguished triangles
$$
\xymatrix{
R\Gamma(X;A_U) \ar[r] \ar[d] & R\Gamma(X;A_{\bar{U}}) \ar[d] \ar[r] & R\Gamma(X;A_{\bar{U}\setminus U}) \ar[d] \ar[r] &\\
R\Gamma(\bar{U}; A_{U}) \ar[r] & R\Gamma(\bar{U}; A_{\bar{U}}) \ar[r] & R\Gamma(\bar{U}; A_{\bar{U}\setminus U}) \ar[r] &
}
$$
Then we obtain a chain of morphisms in cohomology ($j \in \ZZ$)
$$
\xymatrix{
\cdots \ar[r] & H^j(X;A_U) \ar[r] \ar[d] & H^j(X;A_{\bar{U}}) \ar[d] \ar[r] & H^j(X;A_{\bar{U}\setminus U}) \ar[d] \ar[r] & \cdots \\
\cdots \ar[r] & H^j(\bar{U}; A_{U}) \ar[r] & H^j(\bar{U}; A_{\bar{U}}) \ar[r] & H^j(\bar{U}; A_{\bar{U}\setminus U}) \ar[r] & \cdots
}
$$
For $Z=\bar{U}$ or $Z= \bar{U}\setminus U$ we have
    $$R\Gamma(\overline{U};A_Z) \simeq R\Gamma(X;A_Z) \simeq R\Gamma(Z;A_Z)$$
since $A_Z \simeq i_{Z*}i_Z^{-1}A_Z$ (\cite[Proposition 2.3.6 (iv) and (iii)]{ks1}). This implies isomorphisms in cohomology 
$$H^j(\overline{U};A_Z) \simeq H^j(X;A_Z) \simeq H^j(Z;A_Z),$$ for each $j \in \ZZ$.
Therefore by the five lemma we have isomorphisms 
$$H^j(\bar{U};A_U) \simeq H^j(X;A_U)$$ for each $j \in \ZZ.$ 

In the same way, working in $\Dfs$
with the corresponding morphisms of distinguished triangles
 we obtain the isomorphism $$H^j(\bar{U}(\mathbb{S});A_{U(\mathbb{S})}) \simeq H^j(X(\mathbb{S});A_{U(\mathbb{S})})$$ for each $j \in \ZZ.$ 
\qed \\

\begin{lem}\label{nrmk inv Xs to Us}
Let $X$ be an object of $\Df$ and let $U\in \op (X_{\df}).$ If  the morphism in  \eqref{eqn invariance} induces  isomorphisms $H^*(Y;A_Y)\simeq H^*(Y(\mathbb{S});A_{Y(\mathbb{S})})$ for $Y=X$ and $Y= X\setminus U,$ then it also induces  isomorphisms $$H^*(X;A_U)\simeq H^*(X(\mathbb{S});A_{U(\mathbb{S})}).$$
\end{lem}

\pf 
The exact sequence $0 \to A_U \to A_X \to A_{X \setminus U} \to 0$ (\cite[Proposition 2.3.6 (v)]{ks1}) defines a distinguished triangle $A_U \to A_X \to A_{X \setminus U} \to$ (\cite[Proposition 1.7.5]{ks1}). We have the following morphism of distinguished triangles
$$
\xymatrix{
R\Gamma(X;A_U) \ar[r] \ar[d] & R\Gamma(X;A_X) \ar[d] \ar[r] & R\Gamma(X;A_{X \setminus U}) \ar[d] \ar[r] &\\
R\Gamma(X(\mathbb{S});A_{U(\mathbb{S})}) \ar[r] & R\Gamma(X(\mathbb{S});A_{X(\mathbb{S})}) \ar[r] & R\Gamma(X(\mathbb{S});A_{(X \setminus U)(\mathbb{S})}) \ar[r] &
}
$$
where the vertical morphisms are given in \eqref{eqn invariance}.
Then we obtain a chain of morphisms in cohomology ($j \in \ZZ$)
{\tiny
$$
\xymatrix{
\cdots \ar[r] & H^j(X;A_U) \ar[r] \ar[d] & H^j(X;A_X) \ar[d] \ar[r] & H^j(X;A_{X\setminus U}) \ar[d] \ar[r] & \cdots \\
\cdots \ar[r] & H^j(X(\mathbb{S}); A_{U(\mathbb{S})}) \ar[r] & H^j(X(\mathbb{S}); A_{X(\mathbb{S})}) \ar[r] & H^j(X(\mathbb{S}); A_{(X\setminus U)(\mathbb{S})}) \ar[r] & \cdots
}
$$
}
If we set $Z=X\setminus U$, then  we have
    $$R\Gamma(X;A_Z) \simeq R\Gamma(Z;A_Z)$$
    since $A_Z \simeq i_{Z*}i_Z^{-1}A_Z$ (\cite[Proposition 2.3.6 (iv) and (iii)]{ks1}).
In the same way, working in $\Dfs$ and using $X({\mathbb S})\setminus U({\mathbb S})=(X\setminus U)({\mathbb S})$, we have
$$R\Gamma(X(\mathbb{S});A_{Z(\mathbb{S})}) \simeq R\Gamma(Z(\mathbb{S});A_{Z(\mathbb{S})}).$$
Therefore, if  $R\Gamma(Y;A_Y) \simeq R\Gamma(Y(\mathbb{S});A_{Y(\mathbb{S})})$ for $Y=X$ and $Y= X \setminus U,$ then we have isomorphisms in cohomology 
$$H^j(X;A_Y) \simeq H^j(Y;A_Y) \simeq H^j(Y(\mathbb{S});A_{Y(\mathbb{S})}) \simeq H^j(X(\mathbb{S});A_{Y(\mathbb{S})}),$$ $j \in \ZZ$.
Therefore by the five lemma we have isomorphisms 
$$H^j(X;A_U) \simeq H^j(X(\mathbb{S});A_{U(\mathbb{S})})$$
 for each $j \in \ZZ.$ 
\qed \\

We are ready to show our criterion for invariance       of o-minimal sheaf cohomology without supports:

\begin{crt} \label{crt inv without c}
Let $X$ be an object of  $\Df.$ Suppose that $X$ 
has a finite cover by open definable subsets $\{W_i\}_{i=1}^m$ such that each $\bar{W_i}$ is affine.
Suppose that for every (affine)  closed definable subset $Z$ of $X$ with $Z\subseteq \bar{W_i}$ for some $i$ we have an isomorphism
$$H^*(Z;A_{Z})\simeq H^*(Z({\mathbb S}); A_{Z({\mathbb S})}).$$
Then for every  $F \in \mod(A_{X_{\df}})$  we have an isomorphism
$$H^*(X; F)\simeq H^*(X({\mathbb S}); F({\mathbb S})).$$
\end{crt}

\pf
Set $\mathfrak{S}=\{F \in \mod(A_{X_{\df}}):\,   R\Gamma(X;F)\simeq R\Gamma(X(\mathbb{S});F(\mathbb{S}))\}$. We will obtain the result applying Fact \ref{fact  criteria for o-min sheaves}.


The family $\mathfrak{S}$ satisfies (i) and (ii) of Fact \ref{fact  criteria for o-min sheaves}. (i) first is standard: the exact sequence $0 \to F' \to F \to F'' \to 0$ implies the following morphism of distinguished triangles
$$
\xymatrix{
R\Gamma(X;F') \ar[r] \ar[d] & R\Gamma(X;F) \ar[d] \ar[r] & R\Gamma(X;F'') \ar[d] \ar[r] &\\
R\Gamma(X(\mathbb{S});F') \ar[r] & R\Gamma(X(\mathbb{S});F) \ar[r] & R\Gamma(X(\mathbb{S});F'') \ar[r]
&
}
$$
where the vertical morphisms are given in \eqref{eqn invariance}. Therefore, if  $F' \in \mathfrak{S},$ then using the five lemma in the corresponding chain of morphisms in cohomology  we have $F \in \mathfrak{S}$ if and only if $F'' \in \mathfrak{S}$. (ii) is a consequence of the fact that sections commute with filtrant $\Lind$ (\cite[Example 1.1.4 and Proposition 1.2.12]{ep2}).

So we are reduced to proving that the family $\mathfrak{S}$ satisfies (iii)$'$ of  Fact \ref{fact  criteria for o-min sheaves}.  
Namely, we have to show that, for any affine open subset $U$ of $X$ such that $U\subseteq W_i$ for some $i,$ 
we have $$H^*(X;A_U)\simeq H^*(X({\mathbb S}); A_{U({\mathbb S})}).$$ By the isomorphisms  of Lemma \ref{nrmk inv X to U} it is enough to see that we have isomorphisms
 $$H^*(\bar{U};A_U)\simeq H^*(\bar{U}({\mathbb S}); A_{U({\mathbb S})}).$$
 By Lemma \ref{nrmk inv Xs to Us}, we are reduced to proving the isomorphism 
 $$H^*(Y;A_Y) \simeq H^*(Y(\mathbb{S});A_{Y(\mathbb{S})})$$ 
 for $Y=\bar{U}$ and $Y= \bar{U} \setminus U,$ which follows from the fact that, by the hypothesis,  invariance       with constant coefficients holds for every affine closed definable subsets of $X$ contained in some $\bar{W_i}.$ 
\qed \\

\end{subsection}

\begin{subsection}{Criterion for finiteness results}\label{subsection finiteness}
Here we prove our general criterion for our finiteness results, namely  the o-minimal analogue of Wilder's finiteness theorem (\cite[III.10]{i}).\\


\begin{lem}\label{lem wilder tool}
Let $A$ be a  noetherian ring and let $L$ be a finitely generated $A$-module. Let $X$ be an object of $\Df.$ Suppose that $X$ is definably normal. 
Suppose that for every  affine definably compact subset   of $X$ has an affine definably compact neighborhood $B$ such that  $H^q(B;L_B)$ is finitely generated for each $q$.
Then for every   pair $(Z,K)$ of definably compact definable subsets of $X$ such that $K\subseteq \mathring{ Z}$, the restriction map
$$H^q(Z;L_Z)\to H^q(K;L_K)$$
has finitely generated image for each $q$.
\end{lem}

\pf
The proof in by induction on $q$. The result holds for $q<0$ since $H^q(Z;L_Z)=H^q(K;L_K)=0.$ Assume the result holds in degrees $<q.$ Let $\A$ be the collection of all definably compact definable subsets $A$ of $X$ for which there exists a definably compact subset $C$ of $X$ with $A\subseteq \mathring{C}\subseteq C\subseteq \mathring{Z}$ such that the restriction map
$H^q(Z;L_Z)\to H^q(C;L_C)$ has finitely generated image.

\begin{clm}\label{clm wilder}
The collection $\A$ has the following properties:
\begin{enumerate}
\item
If $A$ is an affine definably compact definable subset of $X$ such that $A\subseteq \mathring{ Z},$  then $A\in \A$.

\item
If $A\in \A$ and $R\subseteq A$ is a definably compact subset of $A$, then  $R\in \A$.

\item
If $A\in \A$ and $R\in \A$, then $A\cup R\in \A.$
\end{enumerate}
\end{clm}

\medskip
We obtain (1) by the assumption of affine definably compact subsets of $X$ and,  
(2) is clear. For (3), suppose that $A\in \A$ and $R\in \A$. By definition of $\A$ and the shrinking lemma, there are definably compact subsets $B$ and $C$ such that $A\subseteq \mathring{B}\subseteq B\subseteq  \mathring{C}\subseteq C\subseteq \mathring{Z}$ such that the restriction map $H^q(Z;L_Z)\to H^q(C;L_C)$ has finitely generated image. Similarly, by definition of $\A$ and the shrinking lemma, there are definably compact subsets $S$ and $T$ such that $R\subseteq \mathring{S}\subseteq S\subseteq  \mathring{T}\subseteq T\subseteq \mathring{Z}$ such that the restriction map $H^q(Z;L_Z)\to H^q(T;L_T)$ has finitely generated image. Consider the following commutative diagram constructed from the Mayer-Vietoris sequences (\cite[Chapter II, Section 13 (32) (b)]{b})
{\small
$$
\xymatrix{
                            &  H^q(Z; L_Z) \ar[d] \ar[r] & H^{q}(Z;L_Z)\oplus H^q(Z;L_Z) \ar[d]  &\\
H^{q-1}(C\cap T;L_{C\cap T}) \ar[d] \ar[r] & H^q(C\cup T;L_{C\cup T}) \ar[d] \ar[r] &  H^{q}(C;L_C)\oplus H^q(T;L_T)&\\
H^{q-1}(B\cap S;L_{B\cap S})  \ar[r] & H^q(B\cup S;L_{B\cup S}).   &
}
$$
}

Note that: (i) the middle  horizontal sequence of the diagram  is exact; the first down arrow on the bottom square of the diagram has finitely generated image (by the induction hypothesis); the second down arrow  on the top square of the diagram has finitely generated image (by the hypothesis of (3)). By  the purely algebraic result \cite[III. Lemma 10.3]{i}, we conclude that the restriction map $H^q(Z;L_Z) \to H^q(B\cup S;L_{B\cup S})$ has finitely generated image 
and hence $A\cup R\in \A.$\\

Now let $(Z, K)$ be a pair of  definably compact definable subsets of $X$ such that $K\subseteq \mathring{ Z}.$ Since $X$ is definably normal, by the shrinking lemma, we have $K=K_1\cup \cdots \cup K_r$ where each $K_i$ is an affine  definably compact subset of $X$ such that $K_i\subseteq \mathring{ Z}.$ 
We conclude the proof by induction on $r.$ The case  $r=1$ follows  by Claim \ref{clm wilder} (1) and the inductive step follows by Claim \ref{clm wilder} (3).
\qed \\

Applying  Lemma \ref{lem wilder tool} with $K=Z=X$ we obtain:

\begin{crt}\label{crt wilder tool}
Let $A$ be a  noetherian ring and let $L$ be a finitely generated $A$-module. Let $X$ be an object of $\Df.$ Suppose that $X$ is  definably compact and definably normal.   Suppose that for every  affine closed definable subset $B$  of $X$,  $H^q(B;L)$ is finitely generated for each $q$.
Then $H^q(X;L)$ is finitely generated for each $q$.\\
\end{crt}


\end{subsection}

\end{section}

\begin{section}{Applications to definably compact groups}\label{section applications def grps}
In this section we apply our general invariance      
and finiteness  criteria to obtain, for definably compact definable groups in arbitrary o-minimal structures,  the invariance      
and finiteness results stated in the Introduction. \\

Here we assume that  ${\mathbb M}$ is an arbitrary o-minimal structure and,  as before, we let ${\mathbb S}$ be an elementary extension  of ${\mathbb M}$ or  an o-minimal expansion of ${\mathbb M}$. \\

To proceed we require  the following (\cite[Definition 3.1]{epr}):

\begin{defn}\label{defn gp-int}
{\em
A {\it definable group-interval}  $J=\langle (-b, b), 0, +,  <\rangle $ is an open interval $(-b, b)\subseteq M$, with $-b<b$ in $M\cup \{-\infty, +\infty \},$ together with a binary partial continuous definable operation $+:J^2\to J$ and an element $0\in J$, such that:
\begin{itemize}
\item
 $x+y=y+x$ (when defined), $(x+y)+z=x+(y+z)$ (when defined) and $x<y \Rightarrow x+z<y+z$ (when defined);
\item
for every $x\in J$ with $0<x$, the set $\{y\in J: 0<y \,\,\textrm{and}\,\, x+y \,\, \textrm{is defined}\}$ is an interval of the form $(0, r(x));$
\item
for every $x\in J$ with $0<x$, then $\lim _{z\to 0}(x+z)=x$ and  $\lim _{z\to r(x)^-}(x+z)=b;$
\item
for every $x\in J$ there exists $z\in J$ such that $x+z=0.$
\end{itemize}
The definable group-interval $J$ is {\it unbounded} (resp. {\it bounded}) if the operation $+$ in $J$ is total (resp. not total). The notion of a {\it definable homomorphism} between definable group-intervals is defined in the obvious way.
}
\end{defn}

By the properties above, it follows  that: (i) for each $x\in J$ there is a unique $z\in J$ such that $x+z=0$, called the inverse of $x$ and denoted by $-x$; (ii) for each $x\in J$ we have $-0=0$, $-(-x)=x$ and $0<x$ if and only if $-x<0$; (iii) the maps $J\to J:x\mapsto -x$ and $(-b,0)\to (0,b):x\mapsto -x$  are continuous definable    bijections; (iv) for every $x\in J$ with $x<0$, the set $\{y\in J: y<0 \,\,\textrm{and}\,\, x+y \,\, \textrm{is defined}\}$ is an interval of the form $(-r(x), 0);$ (v) for every $x\in J$ with $x<0$, then $\lim _{z\to 0}(x+z)=x$ and  $\lim _{z\to -r(x)^+}(x+z)=-b;$ (vi) for every $x\in J$ we have $x+0=x$ (both sides are defined and they are equal). \\

By the proof of \cite[Lemma 3.5]{epr} we have:

\begin{fact}\label{fact grp-int grp by 4}
Let $J=\langle (-b, b), 0, +, -,  <\rangle $ is a definable group-interval. Then there exists an injective, continuous definable homomorphism $\tau :J\to J$ given by $\tau (x)=\frac{x}{4}$ such that if  $x, y\in \tau (J)=(-\frac{b}{4}, \frac{b}{4}),$ then $x+y,$  $x-y$ and $\frac{x}{2}$ are defined in $J.$  \\
\end{fact}

From now on we fix a cartesian product  $\bJ =\Pi _{i=1}^mJ_i$  of definable group-intervals $J_i=\langle (-_ib_i, b_i), $ $0_i, +_i, -_i,  <\rangle $. \\

We say that $X$  is a {\it $\bJ $-bounded subset}    if $X\subseteq \Pi _{i=1}^m[-_ic_i, c_i]$ for some $c_i>0_i$ in $J_i$.

\begin{nrmk}\label{nrmk j-bd}
{\em
We will often identify a  $\bJ $-bounded subset $X$ with its image under the cartesian product of the injective homomorphisms given by Fact \ref{fact grp-int grp by 4} and assume that $X\subseteq \Pi _{i=1}^m[-_ic_i, c_i]$ for some $0_i<c_i<\frac{b_i}{4}$ in $J_i.$
}
\end{nrmk}


Let $l\in \{1, \dots , m-1\}$. For a definable subset $X\subseteq  \Pi _{i=1}^lJ_i$, we set $L^l(X)=\{f:X\to   J_{l+1}: f\, \textrm {is definable and continuous}\}$ and $L^l_\infty(X)=L^l(X)\cup\{-_{l+1}b_{l+1}, b_{l+1}\},$ where we regard $-_{l+1}b_{l+1}$ and $b_{l+1}$ as constant functions on $X$.
If $f\in L^l(X)$, we denote by $\Gamma(f)$ the graph of $f$. If $f, g\in L^l_\infty(X)$ with $f(x)<g(x)$ for all $x\in X$, we write $f<g$ and  set $(f,g)_X=\{(x, y)\in X\times J_{l+1}: f(x)<y<g(x)\}$. Then,
\begin{itemize}
 \item
 a \emph{$\bJ$-cell} in $J_1$ is either a singleton subset of $J_1$, or an open interval with endpoints in $J_1\cup \{-_{}b_{}, b_{1}\}$,
 \item
 a \emph{$\bJ$-cell} in $\Pi _{i=1}^{l+1}J_i$ is a set of the form $\Gamma(f)$, for some $f\in L^l(X)$, or $(f,g)_X$, for some $f,g\in L^l_\infty(X)$, $f<g$, where $X$ is a $\bJ$-cell in $\Pi _{i=1}^{l}J_i$.
\end{itemize}
In either case, $X$ is called \emph{the domain} of the defined cell. The {\it dimension} of a $\bJ$-cell in $\Pi _{i=1}^{m}J_i$ is defined as usual (\cite[Chapter 3 (2.3) and Chapter  4 (1.1)]{vdd}).\\

We refer the reader to \cite[Chapter 3 (2.10)]{vdd} for the definition of a \emph{decomposition} of $\bJ .$ A \emph{$\bJ $-decomposition} is then a decomposition $\mathcal{C}$ of $\bJ$ such that each $B\in \mathcal{C}$ is a $\bJ$-cell. The following can be proved similarly to \cite[Chapter 3 (2.11)]{vdd}.

\begin{thm}[$\bJ$-CDT]\label{thm JmCDT}
$\,\,$
\begin{enumerate}
\item
If $A_1, \dots, A_k\subseteq \Pi _{i=1}^mJ_i$ are definable sets, then there is a $\bJ$-decomposition $\mathcal{C}$ that partitions each $A_i$.

\item
If $A\subseteq \Pi _{i=1}^{m}J_i$ is a definable set and $f\in L^{m+1}(A)$, then there is a $\bJ$-decomposition $\mathcal{C}$  that partitions $A$ such that the restriction $f_{|B}$ to each $B\in\mathcal{C}$ with $B\subseteq  A$ is continuous. \\
\end{enumerate}
\end{thm}

Below we will need the following observations. \\

To $\bJ$ there is an associated definable o-minimal structure $\mJ$ such that: (i) the domain of $\mJ$ is the definable set ${\rm dom}(\mJ)=(-_1b_1, b_1)\cup \{c_2\}\cup (-_2b_2, b_2)\cup \ldots \cup \{c_{m}\}\cup (-_mb_m, b_m)$ where the $c_i$'s are new elements, with the obvious induced definable total order; (ii) the $\mJ$-definable subsets are the  subsets $X\subseteq {\rm dom}(\mJ)^k$ such that $X$ is a definable set. \\

By \cite[Fact 4.4]{epr} we have:

\begin{fact}\label{fact Jm def skolem}
The o-minimal structure $\mJ$ has  $\mJ$-definable choice.\\
\end{fact}

Later we will require also the following:

\begin{lem}\label{lem compact is J-bounded}
Let $X\subseteq \Pi _{i=1}^mJ_i$ be a definable subset. Then $X$ is definably compact if and only if $X$ is a closed and $\bJ$-bounded subset.
\end{lem}

\pf
If $X$ is a closed and $\bJ$-bounded subset, then it is a closed and bounded definable set and so, by \cite[Theorem 2.1]{ps}, $X$ is definably compact. 

If $X$ is definably compact, then $X$ is closed and if it were not a $\bJ$-bounded subset, for some $l\in \{1,\ldots , m\},$ the projection of $X$ onto $J_l$ would not be a $J_l$-bounded subset. But then using $\mJ$-definable choice (Fact \ref{fact Jm def skolem}) and arguing as in the proof of \cite[Chapter 6, (1.9)]{vdd} we would contradict the definable compactness of $X.$
\qed \\

The following remark will allow us to work in $\mJ$ instead of in ${\mathbb M}$ when convenient:

\begin{nrmk}\label{nrmk Jm cohom}
{\em
Let $X\subseteq \Pi _{i=1}^mJ_i$ be a definable subset. Then $X$ is a $\mJ$-definable set.
In particular: (i) the o-minimal site of $X$ in ${\mathbb M}$ is the same as the o-minimal site of $X$ in $\mJ$; (ii) the o-minimal cohomology of $X$ computed in ${\mathbb M}$ is the same as the o-minimal cohomology of $X$ computed in $\mJ.$ \\
}
\end{nrmk}

Below we let  $L$ be  an $A$-module.\\

As in the case of o-minimal expansions of ordered groups (\cite[Corollary 3.3]{bf}) we have:

\begin{lem}\label{lem bf acyclic cells}
Let $C$ be a $\bJ$-cell which is a $\bJ$-bounded subset. Then $C$ is  acyclic, i.e. $H^p(C;L_C)=0$ for $p>0$ and $H^0(C; L_C)=L$.
\end{lem}

\pf
This is obtained in exactly the same way as \cite[Corollary 3.3]{bf}. Indeed, since $C$ is a $\bJ$-bounded subset, by Remark \ref{nrmk j-bd}, we can apply the group-interval operations $x+_iy$, $x-_iy$ and $\frac{x}{2}$ in each coordinate of $ \Pi _{i=1}^mJ_i$ just like in the proof of \cite[Corollary 3.3]{bf} obtaining:

\begin{clm}\label{clm bd int}
If $I$ is a definably connected $J_1$-bounded subset, then $I$ is definably contractible to a point in $J_1$.
\end{clm}

\begin{clm}\label{clm bd cell retract}
If $C$ is a $\bJ$-cell which is a $\bJ$-bounded subset, then there is a definable deformation retract of $C$ to a $\bJ$-cell which is a $\bJ$-bounded subset of strictly lower dimension.
\end{clm}
\noindent
See \cite[Lemmas 3.1 and 3.2]{bf}.

By Claim \ref{clm bd cell retract} and induction on the dimension of $C$,  $C$ definably contractible to a point in $ \Pi _{i=1}^mJ_i.$ Note also that by construction the  domain of the definable deformation retraction of Claim \ref{clm bd cell retract} is a definable subset of $\Pi _{i=1}^mJ_i.$ Therefore,  by Fact \ref{fact Jm def skolem} and Remark \ref{nrmk Jm cohom}, we have the homotopy axiom for o-minimal cohomology (\cite{ejp}) for definable homotopies whose domains are definable subsets of $\Pi _{i=1}^mJ_i.$ So $H^p(C;L_C)$ is the same as the o-minimal cohomology of a point and we apply the dimension axiom for o-minimal cohomology to conclude. \qed \\

We also have the analogue of \cite[Lemma 7.1]{bf}:

\begin{lem}\label{lem bf cells}
Let $C$ be a $\bJ$-cell which is a $\bJ$-bounded subset and  of dimension $r$. There is a definable family $\{C_{t_1,\ldots , t_m} : 0_i<t _i<\frac{b_i}{4}, \,\, i=1, \ldots , m\}$ of closed and $\bJ$-bounded subsets $C_{t_1,\ldots , t_m} \subset C$ such that:
\begin{enumerate}
\item
$C = \bigcup_{t_1,\ldots, t_m}C_{t_1,\ldots , t_m}.$
\item If $0_i < t'_i < t_i$ for all $i=1, \ldots , m$, then $C_{t_1,\ldots , t_m}
\subset C_{t'_1,\ldots , t'_m}$ and this inclusion induces an isomorphism 
$$H^{p}(C \backslash C_{t_1,\ldots , t_m}; L_C) \simeq  H^{p}(C \backslash C_{t'_1,\ldots , t'_m}; L_C).$$
\item
The o-minimal cohomology of $C \backslash C_{t_1,\ldots , t_m}$ is given by
\begin{equation*}
H^{p}(C \backslash C_{t_1,\ldots , t_m}; L_C) =
\begin{cases}
L ^{1+\chi _{1}(r)}\qquad \textmd{if} \qquad p\in \{0, r-1\}\\
\\
\,\,\,\,\,\,\,\,\,\,\,\,\,\,\,\,\,\,\,\,\,\,\,\,\,\,\,\,\,\,\,
\\
0\qquad \,\,\,\,\,\,\, \,\,\,\,\,\,\,\,\,\,\, \, \textmd{if} \qquad p\notin \{0, r-1\}
\end{cases}
\end{equation*}
where $\chi _1:{\mathbb Z}\to \{0,1\}$ is the characteristic function of the subset $\{1\}.$
\end{enumerate}
\end{lem}

\pf
By Remark \ref{nrmk j-bd}, we assume that $C\subseteq \Pi _{i=1}^m[-_ic_i, c_i]$ for some $0_i<c_i<\frac{b_i}{4}$ in $J_i$ and the group-interval operations $x+_iy$, $x-_iy$ and $\frac{x}{2}$ are all defined in each coordinate of $ \Pi _{i=1}^mJ_i.$

We define the definable family $\{C_{t_1,\ldots , t_m} : 0_i<t _i<\frac{b_i}{4}, \,\, i=1, \ldots , m\}$ by induction on $l\in \{1, \ldots , m-1\}$ in the following way.
\begin{enumerate}
\item
If $l=1$ and $C$ is a singleton in $J_1$, we define $C_{t_1} = C.$
\item
If $l= 1$ and $C = (d, e)\subseteq J_1,$ then $C_{t_1} =[d +_1 \gamma ^1_{t_1}, e -_1\gamma ^1_{t_1}]$ where $\gamma ^1_ {t_1} = \min\{|\frac{d-_1e}{2}|_1, t_1\},$ (in this way $C_{t_1}$ is non empty).
\item
If $l> 1$ and $C = \Gamma(f),$ where $f\in L^l(B)$ is a continuous definable map and $B$ is $\bJ$-cell in $\Pi _{i=0}^lJ_i$ which is a $\bJ$-bounded subset. By induction $B_{t_1, \ldots, t_l}$ is defined. We put $C_{t_1, \ldots, t_l, t_{l+1}}=\Gamma (f_{|B_{t_1,\ldots , t_l}}).$

\item
If $l > 1$ and $C = (f, g)_{B},$ where $f, g\in L^l(B)$ are continuous definable maps,  $B$ is $\bJ$-cell in $\Pi _{i=0}^lJ_i$ which is a $\bJ$-bounded subset and  $f<g$. By induction $B_{t_1,\ldots, t_l}$ is defined. We put $C_{t_1,\dots , t_l, t_{l+1}} = [f +_{l+1}\gamma ^{l+1}_{t_{l+1}}, g -_{l+1} \gamma ^{l+1}_{t_{l+1}}]_{B_{t_1,\ldots , t_l}},$ where $\gamma ^{l+1}_{t_{l+1}} := \min(|\frac{f-_{l+1} g}{2}|_{l+1}, t_{l+1}).$
\end{enumerate}

We observe that from this construction we obtain:

\begin{clm}\label{clm bf cover cells}
For $t_1,\ldots , t_m $ as above there is a covering ${\mathcal{U}}_{C} = \{U_i : i \in I\}$ of $C \backslash C_{t_1,\ldots , t_m}$ by relatively open $\bJ$-bounded subset such that:
\begin{enumerate}
\item
The index set $I$ is the family of the closed faces of an $r$-dimensional cube. (So $|I| = 2r$).
\item
If $E \subset I,$ then $U_{E} := \bigcap_{i \in E}U_{i}$ is either empty or a $\bJ$-cell. (So in particular $H^{p}(U_{E}; L_C) = 0$ for $p > 0$ and, if $U_{E} \neq \emptyset, H^{0}(U_{E}; L_C) = L.$)
\item
For $E \subset I,$ $U_{E} \neq \emptyset$ iff the faces of the cubes belonging to $E$ have a non-empty intersections.
\end{enumerate}
So the nerve of ${\mathcal U}_C$ is isomorphic to the nerve of a covering of an $r$-cube by its closed faces.
\end{clm}

\pf
To show that there is a covering satisfying the properties above, we define ${\mathcal{U}}_C$ by induction on $l\in \{1, \ldots , m-1\}$. We distinguish four cases
according to definition of the $C_{t_1,\ldots , t_m}.$

\begin{enumerate}
\item
If $l= 1$ and $C$ is a singleton in $J_1$, then $\mathcal{U}_{C}$ is the covering consisting of one open set (given by the whole space $C$).
\item
If $l= 1$ and $C = (d, e)\subseteq J_1,$ then $C \backslash C_{t_1}$ is the union of the two open subsets $(d, d +_1 \gamma ^1_{t_1})$ and $(e -_1 \gamma ^1_{t_1}, e),$ and we define ${\mathcal{U}}_C$ as the covering consisting of these two sets.
\item
If $l > 1$ and $C = \Gamma(f),$ where $f\in L^l(B)$ is a continuous definable map and $B$ is $\bJ$-cell in $\Pi _{i=0}^lJ_i$ which is a $\bJ$-bounded subset. By definition $C_{t_1, \ldots, t_l, t_{l+1}}=\Gamma (f_{|B_{t_1,\ldots , t_l}}).$ By induction we have a
covering $\mathcal{V}_B$ of $B \backslash B_{t_1, \ldots , t_l}$ with the stated
properties, and we define ${\mathcal{U}}_C$ to be a covering of $C \backslash C_{t_1,\ldots , t_l, t_{l+1}}$ induced by the natural homeomorphism between the graph of $f$ and its domain.
\item
If $l > 1$ and $C = (f, g)_{B}$ where $f, g\in L^l(B)$ are continuous definable maps,  $B$ is  $\bJ$-cell in $\Pi _{i=0}^lJ_i$ which is a $\bJ$-bounded subset and  $f<g$. By definition $C_{t_1,\dots , t_l, t_{l+1}} = [f +_{l+1}\gamma ^{l+1}_{t_{l+1}}, g -_{l+1} \gamma ^{l+1}_{t_{l+1}}]_{B_{t_1,\ldots , t_l}}.$ By induction we have that $B \backslash B_{t_1, \ldots , t_l}$ has a covering
${\mathcal V}_B = \{V_{j} : j \in J\}$ with the stated properties, where $J$ is the set of closed faces of the cube $[0, 1]^{r-1}.$
Define a covering ${\mathcal{U}}_C = \{U_i : i \in I\}$ of $C \backslash C_{t_1,\dots , t_l, t_{l+1}}$ as follows. As index set $I$ we take the closed faces of the cube $[0, 1]^r.$ Thus $|I|=|J|+2,$ with the two extra faces  corresponding to the ``top" and ``bottom" face of $[0, 1]^r.$ We
associate to the top face the open set $(g -_{l+1} \gamma ^{l+1}_{t_{l+1}}, g)_{B_{t_1, \ldots , t_l}}
\subset C \backslash C_{t_1,\dots , t_l, t_{l+1}}$ and the bottom face the open set $(f, f
+_{l+1} \gamma ^{l+1}_{t_{l+1}})_{B_{t_1, \ldots , t_l}} \subset C \backslash C_{t_1,\dots , t_l, t_{l+1}}.$ The other open sets of the covering are the preimages of the sets $V_j$ under the restriction of the projection $\Pi _{i=1}^{l+1}J_i \rightarrow \Pi _{i=1}^lJ_i.$ This defines a covering of $C \backslash C_{t_1,\dots , t_l, t_{l+1}}$ with the stated properties.
\end{enumerate}
\qed

Property (1) of the lemma is clear. By (the proof of) Claim \ref{clm bf cover cells} there are open covers ${\mathcal U}'_C$ of $C \backslash C_{t'_1,\ldots , t'_m}$ and ${\mathcal U}_C$ of $C \backslash C_{t_1,\ldots , t_m}$ satisfying the assumptions of \cite[Lemma 5.5]{bf}. Hence property (2) of the lemma holds. Finally, if $r>1$, then property (3) follows from  Claim \ref{clm bf cover cells}  and \cite[Corollary 5.2]{bf}. On the other hand, if $r=1$, then $C\setminus C_{t_1,\ldots , t_m}$ is by construction a disjoint union $D\sqcup E$ of two $\bJ$-cells which are $\bJ$-bounded subsets and  of dimension $r=1$. Therefore, in this case, the result follows from Lemma \ref{lem bf acyclic cells}, since $H^{*}(C \backslash C_{t_1,\ldots , t_m}; L_C) \simeq H^*(D;L_D)\oplus H^*(E;L_E).$
\qed  \\

From Lemma \ref{lem bf cells} and computations in o-minimal cohomology we obtain just like in \cite[Lemma 7.2 and Corollary 7.3]{bf}:

\begin{lem}\label{lem bf x bd}
Let $X$ be a  definable  $\bJ$-bounded subset  and $C\subseteq X$  a $\bJ$-cell of maximal dimension. Then for every $t_1,\ldots , t_m $ and $t'_1,\ldots , t'_m $ with $t_i'<t_i$ for all $i=1,\ldots , m$ as  above we have isomorphisms induced by inclusions:

\begin{enumerate}
\item
$H^{*}(X \backslash C_{t_1,\ldots , t_m}; L_X) \simeq H^{*}(X \backslash C_{t'_1,\ldots , t'_m}; L_X);$
\item
$H^{*}(X \backslash C_{t_1,\ldots , t_m}; L_X) \simeq H^{*}(X \backslash C; L_X)$ assuming also that $X$ is closed.
\end{enumerate}
\end{lem}



\begin{nrmk}\label{nrmk Uct}
{\em
Let $X$ be a  definable  $\bJ$-bounded subset  and $C\subseteq X$  a $\bJ$-cell. Assume that $C\subseteq \Pi _{i=1}^m[-_ic_i, c_i]$ for some $0_i<c_i<\frac{b_i}{4}$ in $J_i.$ Then there is a point $p_C\in C$ such that for all $t_1,\ldots , t_m $ as  above, if   $c_i<t_i$ for all $i=1,\ldots , m,$  then  $C_{t_1,\ldots , t_m}=\{p_C\}.$ In particular, we have
$$H^*(\bar{C}\setminus \{p_C\};L_X)\simeq H^*(\bar{C}\setminus C ; L_X)$$
even if $\bar{C}$ is in general non-acyclic (\cite[Theorem  4.1]{bf}).
}
\end{nrmk}

From Lemma \ref{lem bf x bd} and computations in o-minimal cohomology we obtain just like in \cite[Theorem 8.1]{bf}:

\begin{lem}\label{lem bf inv gp-int}
If $X$ is a  closed definable $\bJ$-bounded subset, then we have an isomorphism
$$H^*(X;A_{X})\simeq H^*(X({\mathbb S}); A_{X({\mathbb S})}).$$ \\
\end{lem}



Also  we obtain just like in \cite[Theorem 7.4]{bf}:

\begin{lem}\label{lem fg coho def groups}
Let $A$ be a  noetherian ring and let $L$ be a  finitely generated  $A$-module. If $X$ is a closed definable $\bJ$-bounded subset, then $H^p(X;L_X)$ is  finitely generated for each $p.$ \\ 
\end{lem}

We can now go back to the proof of  Theorems \ref{thm inv def groups} and  \ref{thm omin gp wilder} in the Introduction. Below we assume the reader familiarity with the basic theory of definable groups (\cite{ot} and
\cite{p1}). \\

We will require the following result (\cite[Theorem 3]{epr}):

\begin{fact}\label{fact grp epr}
If $G$ is a definable group, then there is a definable injection $\theta: G\to \Pi _{i=1}^mJ_i$, where each $J_i\subseteq M$ is a definable group-interval. \\
\end{fact}

It follows  that:

\begin{lem}\label{lem G is j-man}
If $G$ is a definable group, then  there is a cartesian product $\bJ=\Pi _{i=1}^mJ_i$ of  definable group-intervals  such that $G$ has definable charts $\{(U_i,\phi _i)\}_{i=1}^k$ with $\phi _i(U_i)\subseteq \Pi _{i=1}^mJ_i$  for each $i.$
\end{lem}

\pf
By Fact \ref{fact grp epr}, there is a  cartesian product  $\bJ=\Pi _{i=1}^mJ_i$ of  definable group-intervals such that $G$ is definably isomorphic to a definable group $H\subseteq \bJ.$ By \cite{p1} definable isomorphisms of definable groups are definable homeomorphisms when each definable group is equipped with its definable manifold structure. So we may assume that $G\subseteq \Pi _{i=1}^mJ_i.$ 

By the construction of the definable manifold  structure of $G$ (\cite{p1}), $G$ is a definable space whose definable charts $\{(U_i, \phi _i)\}_{i=1}^k$ are such that each $U_i\subseteq  \Pi_{i=1}^mJ_i$ is a definable subset.  In fact, each $U_i$ is a $\bJ$-cell in $G\subseteq \Pi _{i=1}^mJ_i$ of dimension $n=\dim G$ or $U_i$ is a translate $g_iU_{i'}$ in $G$ of a $\bJ$-cell $U_{i'}$ in $G\subseteq \Pi _{i=1}^mJ_i$ of dimension $n$. In the first case, $\phi _i$ is the restriction of a projection from $\Pi _{i=1}^mJ_i$ onto some $n<m$ coordinates. In the second case $\phi _i$ is the composition of the  translation $x\mapsto g_ix$ in $G$ and the restriction of the projection $\phi _{i'}$ as above. For the fact that the restriction of a projection as above is a definable homeomorphism compare with \cite[Chapter 3, (2.7)]{vdd}.
\qed \\

In the paper \cite{emp} we proved (following a suggestion of the referee of the present paper) that in o-minimal structures with definable Skolem functions, definably compact (Hausdorff) definable spaces are definably normal. For our applications we will need that definably compact definable groups are definably normal and this was previously proved in \cite[Corollary 2.3]{et} and it now follows also from Lemma \ref{lem G is j-man}, Fact \ref{fact Jm def skolem} and the quoted result from \cite{emp}. \\

\begin{thm}\label{thm main J inv and finite}
Let $X$ be a (Hausdorff) definably compact definable space. Suppose that there is a cartesian product $\bJ=\Pi _{i=1}^mJ_i$ of  definable group-intervals  and $X$ has definable charts $\{(U_i,\phi _i)\}_{i=1}^k$ such that $\phi _i(U_i)\subseteq \Pi _{i=1}^mJ_i$  for each $i.$ 
\begin{enumerate}
\item
If $F$ is a sheaf  on the o-minimal site  on $X$, then
we have
$$H^*(X;F)\simeq H^*(X({\mathbb S});F({\mathbb S})).$$
\item
If $A$ is a  noetherian ring and $L$ a  finitely generated  $A$-module, then we have that $H^p(X;L_X)$ is  finitely generated for each $p.$ 
\end{enumerate}
\end{thm}

\pf
Since $X$  is definably normal (as explained above), by the shrinking lemma, $X$ has  a finite cover by  open definable subsets $\{W_i\}_{i=1}^k$  such that $\bar{W_i}\subseteq U_i.$ In particular, each $\bar{W_i}$ is affine. Since $X$ is definably compact, each $\bar{W_i}$ is definably compact. Therefore, each $\phi _i(\bar{W_i})\subseteq  \Pi_{i=1}^mJ_i$ is a definably compact subset. Hence, by Lemma \ref{lem compact is J-bounded}, each $\phi _i(\bar{W_i})$ is a closed and $\bJ$-bounded subset.  In particular,  every (affine)  closed definable subset $Z$ of $X$ with $Z\subseteq \bar{W_i}$ for some $i,$ we have that $\phi _i(Z)\subseteq \phi _i(\bar{W_i})$ is a  closed and $\bJ$-bounded subset. 

Now (1) follows from Lemma \ref{lem bf inv gp-int} and  Criterion \ref{crt inv without c} and (2) follows from Lemma \ref{lem fg coho def groups} and Criterion \ref{crt wilder tool}.
\qed\\

\medskip
\noindent
\textbf{Proofs of Theorems \ref{thm inv def groups} and \ref{thm omin gp wilder}:}
Let $G$ be a definably compact definable group. We have to show respectively that: 
\begin{enumerate}
\item
If $F$ is a sheaf  on the o-minimal site  on $G$, then
we have
$$H^*(G;F)\simeq H^*(G({\mathbb S});F({\mathbb S})).$$
\item
If $A$ is a  noetherian ring and $L$ a  finitely generated  $A$-module, then we have that $H^p(G;L_G)$ is  finitely generated for each $p.$ 
\end{enumerate}
But these follow now at once from Lemma \ref{lem G is j-man} and Theorem \ref{thm main J inv and finite}.
\qed\\

As special case of Theorem \ref{thm main J inv and finite} we obtain the following generalization of the invariance and finiteness results for o-minimal cohomology, without supports, with constant coefficients, on closed and bounded  definable sets in o-minimal expansions of ordered groups (\cite[Theorems 7.4 and  8.1]{bf}):


\begin{cor}\label{cor inv and finite cp exp groups}
Suppose that ${\mathbb M}$ is an o-minimal expansion of an ordered group.  Let $X$ be a (Hausdorff) definably compact,  definable space. 
\begin{enumerate}
\item
If $F$ is a sheaf  on the o-minimal site  on $X,$ then
we have
$$H^*(X;F)\simeq H^*(X({\mathbb S});F({\mathbb S})).$$
\item
If $A$ is a  noetherian ring and $L$ a  finitely generated  $A$-module, then we have that $H^p(X;L_X)$ is  finitely generated for each $p.$ 
\end{enumerate}
\end{cor}


In o-minimal expansions of real closed fields, (1) of the later can be improved due to the following (\cite[Corollary 1.3]{ew}):

\begin{fact}\label{fact ew inv}
If ${\mathbb M}$ is an o-minimal expansion of a real closed field, then for every  definable subset $X\subseteq M^n$ we have an isomorphism
$$H^*(X;A_{X})\simeq H^*(X({\mathbb S}); A_{X({\mathbb S})}).$$
\end{fact}

This fact can be extended to non affine definable spaces and non constant coefficients in the following way:

\begin{cor}\label{cor inv field without c}
Suppose that  ${\mathbb M}$ is an o-minimal expansion of a real closed field. Let $X$ be a regular definable space. Let $F$ be a sheaf  on the o-minimal site  on $X$. Then
we have 
$$H^*(X;F)\simeq H^*(X({\mathbb S}); F({\mathbb S})).$$
\end{cor}

\pf
From \cite[Chapter 10, (1.8) and Chapter 6, (3.5)]{vdd} $X$ is a definably normal definable space. We obtain the result applying Criterion \ref{crt inv without c} using  Fact \ref{fact ew inv}.
\qed\\

\end{section}

\begin{section}{Concluding remarks}\label{section conclusion}
Here we assume that  ${\mathbb M}$ is an arbitrary o-minimal structure. As before, we let ${\mathbb S}$ be an elementary extension  of ${\mathbb M}$ or  an o-minimal expansion of ${\mathbb M}$. We will show the  invariance   results for  o-minimal sheaf cohomology with definably compact supports mentioned in the Introduction. \\

We will need the fact,  mentioned before, and proved in the paper \cite{emp}, that  in o-minimal structures with definable Skolem functions, definably compact (Hausdorff) definable spaces are definably normal. In particular, since existence of definable Skolem functions and Hausdorff are invariant under the functor $\Df \to \Df({\mathbb S})$, definably normal is also invariant. \\


We say that a definable space $X$ is {\it definably completable} if there exists a (Hausdorff) definably compact space $P$ together with a definable open immersion $\iota: X\hookrightarrow P$ (i.e. $\iota (X)$ is open in $P$ and $\iota :X\to \iota (X)$ is a definable homeomorphism)  with  $\iota (X)$ dense in $P$. Such $\iota : X\hookrightarrow P$ is called a {\it definable completion} of $X$. \\

Note that the definition of definable completion appears in the case of o-minimal expansions of real closed fields in the affine case  in \cite[Chapter 10, (2.5)]{vdd} but without the requirement that $\iota(X)$ is open in $P.$ The reason is that, as pointed out in \cite{Fo}, since in that context $P$ is definably normal (\cite[Chapter 6, (3.5)]{vdd}), if $X$ is  locally definably compact (i.e. every point in $X$ has a definably compact neighborhood), then $\iota (X)$ is open in $P.$ By what we said above this is the case here also.\\

\begin{lem} \label{lem inv with c}
Suppose that ${\mathbb M}$ has definable Skolem functions. Let $\iota :X\to P$ be a definable completion of a definable space $X.$ Let $F$ be a sheaf on the o-minimal site on $X$. Then  we have isomorphisms
$$H_c^*(X;F)\simeq H^*(P; \iota _!F)\,\,\,\textrm{and}\,\,\,H_c^*(X({\mathbb S});F({\mathbb S}))\simeq H^*(P({\mathbb S}); (\iota ^{{\mathbb S}})_!F({\mathbb S})).$$
Moreover, $(\iota _!F)({\mathbb S}) \simeq (\iota ^{{\mathbb S}})_!F({\mathbb S})$ (applying tilde we are in the case of diagram \eqref{loc closed base change} on page \pageref{loc closed base change}).
\end{lem}

\pf
Since $P$ is also definably normal (\cite{emp}), by Remark \ref{nrmk def comp2} every definably compact definable subset of $P$ is closed and so  definably normal. Thus  $c,$ the family of definably compact definable subsets of $X,$ is a definably normal  family of supports on $X$. Therefore, by \cite[Corollary 3.9]{ep1} we have $H_c^*(X;F)\simeq H^*(P; \iota _!F).$

Using  the invariance  of definable open immersion and definably compact (Remarks \ref{nrmk def normal inv of compact1} and  \ref{nrmk def normal inv of compact2}) we see that  there is an open ${\mathbb S}$-definable  immersion $\iota ^{{\mathbb S}}: X({\mathbb S})\to P({\mathbb S})$ with $P({\mathbb S})$ an ${\mathbb S}$-definably compact object of $\Dfs$. Since $P({\mathbb S})$ is also ${\mathbb S}$-definably normal (\cite{emp}), by Remark \ref{nrmk def comp2}  in ${\mathbb S}$ every  ${\mathbb S}$-definably compact  ${\mathbb S}$-definable subset of $P({\mathbb S})$ is closed and so  ${\mathbb S}$-definably normal. Thus $c,$ the family of ${\mathbb S}$-definably compact ${\mathbb S}$-definable subsets of $X({\mathbb S}),$ is an  ${\mathbb S}$-definably normal  family of supports on $X({\mathbb S})$. Therefore, by \cite[Corollary 3.9]{ep1}  in ${\mathbb S}$ we have $H_c^*(X({\mathbb S});F({\mathbb S}))\simeq H^*(P({\mathbb S}); (\iota ^{{\mathbb S}})_!F({\mathbb S})).$ 


\qed \\


By Theorem \ref{thm main J inv and finite} and Lemma \ref{lem inv with c} we get:

\begin{cor}\label{cor inv  exp groups with c}
Let $X$ be a  definable space with a definable completion $\iota :X\to P.$ Suppose that there is a cartesian product $\bJ=\Pi _{i=1}^mJ_i$ of  definable group-intervals  and $P$ has definable charts $\{(U_i,\phi _i)\}_{i=1}^k$ such that $\phi _i(U_i)\subseteq \Pi _{i=1}^mJ_i$  for each $i.$  If $F$ is a sheaf  on the o-minimal site  on $X$, then
we have
$$H_c^*(X;F)\simeq H_c^*(X({\mathbb S});F({\mathbb S})).$$\\
\end{cor}






By \cite[Chapter 10, (1.8),  Chapter 10, (2.5)]{vdd}  (and the observation about the definition of definable completable) we have:

\begin{fact}\label{fact comp in fields}
If ${\mathbb M}$ is an o-minimal expansion of a real closed field, then every regular,  locally definably compact definable space  is definably completable by an affine definable space.
\end{fact}

By Facts \ref{fact comp in fields} and Corollary \ref{cor inv  exp groups with c} we also have:

\begin{cor}\label{cor inv fields}
Suppose that ${\mathbb M}$ is an o-minimal expansion of an ordered field. Let $X$ be a   regular,  locally definably compact definable space. Let $F$ be a sheaf  on the o-minimal site on $X$. Then
we have
$$H_c^*(X;F)\simeq H_c^*(X({\mathbb S});F({\mathbb S})).$$
\end{cor}

\end{section}

\end{document}